\documentclass[12pt]{amsart}
\usepackage{amsmath}
\usepackage{amsfonts}
\usepackage{amsthm}
\usepackage{amssymb}
\usepackage{amscd}
\usepackage[all]{xy}
\usepackage{enumerate}


\theoremstyle{plain}
\newtheorem{thm}{Theorem}[subsection]

\newtheorem{prop}[thm]{Proposition}

\newtheorem{cor}[thm]{Corollary}
\newtheorem{lem}[thm]{Lemma}

\theoremstyle{definition}
\newtheorem{defn}[thm]{Definition}

\newtheorem{conv}[thm]{Convention}

\newtheorem*{ackn}{Acknowledgement}

\theoremstyle{remark}
\newtheorem*{rem}{Remark}

\newcommand{\sC}{\mathcal{C}}
\newcommand{\sD}{\mathcal{D}}

\newcommand{\sF}{\mathcal{F}}
\newcommand{\sG}{\mathcal{G}}
\newcommand{\sH}{\mathcal{H}}
\newcommand{\sI}{\mathcal{I}}
\newcommand{\sJ}{\mathcal{J}}
\newcommand{\sL}{\mathcal{L}}
\newcommand{\sN}{\mathcal{N}}
\newcommand{\sM}{\mathcal{M}}
\newcommand{\sO}{\mathcal{O}}

\newcommand{\sS}{\mathcal{S}}

\newcommand{\tA}{{\widetilde{A}}}

\newcommand{\mF}{\mathbb{F}}

\newcommand{\mP}{\mathbb{P}}

\newcommand{\mZ}{\mathbb{Z}}

\newcommand{\Ima}{\mathrm{Im}\,}

\newcommand{\ap}{\mathrm{ap}}

\newcommand{\Aut}{\mathrm{Aut}\,}
\newcommand{\Hilb}{\mathrm{Hilb}}

\newcommand{\Pic}{\mathrm{Pic}\,}
\newcommand{\PGL}{\mathrm{PGL}}

\newcommand{\codim}{\mathrm{codim}\,}

\newcommand{\VSP}{\mathrm{VSP}\,}

\numberwithin{equation}{section}

\title{Spin curves and Scorza quartics} 

\author{Hiromichi Takagi}

\address{Graduate School of Mathematical Sciences \\
the University of Tokyo\\
Tokyo, 153-8914, Japan\\
\texttt{takagi@ms.u-tokyo.ac.jp}}

\author{Francesco Zucconi}

\address{D.I.M.I. \\
the University of Udine\\
Udine, 33100 Italy\\
\texttt{Francesco.Zucconi@dimi.uniud.it}}

\dedicatory{Dedicated to Professor Miles Reid
on the occasion of his 60-th birthday}

\begin{document}

\maketitle

\markboth{Takagi and Zucconi}{Spin curves and Scorza quartics}

\begin{abstract}
In the paper \cite{TZ},
we construct new subvarieties in the varieties of power sums for
certain quartic hypersurfaces.
In this paper, we show that
these quartics coincide with
the Scorza quartics of general
pairs of trigonal curves and ineffective theta characteristics.
Among other applications, we give 
an affirmative answer to the conjecture of Dolgachev and Kanev on 
the existence of the Scorza quartics for any general pairs of
curves and ineffective theta characteristics. 
We also give descriptions of the moduli spaces of
trigonal even spin curves.
\end{abstract}

\tableofcontents

\section{Introduction}
We give a slightly long introduction 
for general readers to be able to learn the essential of 
our new results by reading this introduction only.
Besides we include detailed explanations
of important notions for our paper.
Especially in \ref{subsection:IntroSc},
we give a detailed explanation of the Scorza quartic,
and, in \ref{subsection:Special},
we review our results of \cite{TZ}
needed in this paper.
\subsection{Even spin curves}
\label{section:recover}~

A {\em{theta characteristic}} on a smooth curve $\Gamma$ of genus $g$
is an element $\theta\in {\rm{Pic}}\,\Gamma$ such that $2\theta$ 
is the class of the canonical sheaf $\omega_{\Gamma}$. 
A couple $(\Gamma,\theta)$ is called
a {\it{spin curve}}. 

The study of spin curves started with Riemann
himself 
and since then  a vast 
literature
is devoted to 
the understanding of these objects
(see \cite{Mum}, \cite{topic}).

There are $2^{2g}$ different kinds of spin curve structures
for every smooth curve $\Gamma$ and they are partitioned into two classes according to the parity of $h^{0}(\Gamma,\theta)$.
We say a theta characteristic $\theta$ is {\em{even}} or {\em{odd}}
if $h^{0}(\Gamma,\theta)$ is even or odd respectively.
Correspondingly we speak of
{\em{even or odd spin curves}}.
There exists the moduli space $\sS_{g}$ which parameterizes 
smooth spin curves $(\Gamma,\theta)$
and by the forgetful map $\sS_{g}\rightarrow \sM_{g}$,
where $\sM_{g}$ is the moduli space of curves of genus $g$,
we see that $\sS_{g}$ is a disjoint union of two irreducible
components $\sS^{+}_{g}$ and $\sS^{-}_{g}$ of relative degrees
$2^{g-1}(2^{g}+1)$ and $2^{g-1}(2^{g}-1)$ corresponding to even and
odd theta characteristics respectively.  

It is known that $h^{0}(\Gamma,\theta)=0$ 
for a general pair 
$(\Gamma,\theta)\in \sS^{+}_{g}$ and this lack of sections is a
difficulty to the study of these theta characteristics, which are
called {\it{ineffective theta characteristics}}. On the other hand,
given such an ineffective $\theta$, 
it holds that $h^{0}(\Gamma,\theta+a)=1$ for every $a\in\Gamma$
by the Riemann-Roch theorem, hence $\theta$ gives a
correspondence $I_{\theta}\subset\Gamma\times\Gamma$ such that 
$(a,b)\in I_{\theta}$ if and only if 
$b$ is in the support of the unique member of $|\theta+a|$. This correspondence, called
{\it{the Scorza correspondence}}, is the basis 
for our study in this paper (see \ref{subsection:introtheta}), and
also for two known important applications of even spin curves,
which we now explain. 

The first one is the proof of rationality of $S^{+}_{3}$.
We follow the explanation in \cite[\S 6, 7]{DK}
(see also \cite[\S 3]{Schr}).
Let $V$ be a $3$-dimensional vector space
and $\check{V}$ its dual.
For a homogeneous form $G\in S^m \check{V}$ of degree $m$ on $V$,
we define the (first) polar $P_a(G)$ of $G$ at $a\in \mP(V)$ by
$P_a(G):=
\frac{1}{m}\sum a_i \frac{\partial G}{\partial x_i}$,
where $a_i$ and $x_i$ are coordinates of $a$, and on $V$, respectively.
Let $F\in S^{4}\check{V}$ be a general ternary quartic form on $V$. 
Then the closure of the loci in $\mP(V)=\mP^2$ at a point of 
which the first polar of $F$ is a Fermat cubic is again a smooth
quartic curve,
which is
denoted by $S(F)$ and is called the Clebsch covariant quartic of $F$.
By taking the second polars of
$S(F)$, we have the following correspondence:

\begin{equation}\label{rangouno}
T(F):=\{(a,b)\in S(F)\times S(F)\mid {\rm{rank}}P_{a,b}(F)\leq 1\},
\end{equation}
\noindent which is equal to the correspondence $I_{\theta}$ defined by a unique theta characteristic $\theta$. 
Actually, this is also equal to
$\{(a,b)\in \mP^2 \times \mP^2\mid {\rm{rank}}P_{a,b}(F)\leq 1\}.$

So we have the map
${\rm{Sc}}\colon 
\sM^{0}_{3}\rightarrow  S^{+}_{3}$
such that
${\rm{Sc}}\colon \{F=0\}\mapsto (S(F),\theta)$ 
defined over the open set
$\sM^{0}_{3}\subset \sM_{3}$ where $S(F)$ is nonsingular. 
This association map was discovered by Scorza and is called
the {\em{Scorza map}}.
It turns out to be an injective birational map (\cite[Theorem 7.8]{DK}).
Now we can conclude that $ S^{+}_{3}$ is rational 
since $\sM_{3}$ is rational \cite{Kat} (see also \cite{Boh}).
Nowadays the curve $F$ corresponding to a couple $(S(F),\theta)$ is
called the {\it{Scorza quartic}} of $(S(F),\theta)$.

The second application is Mukai's description of a 
famous prime Fano threefold of genus $12$,
which is, by definition, 
a smooth projective threefold $X$ such that $-K_{X}$ is ample, 
the class of $-K_{X}$ generates $\Pic X$, 
and such that 
the genus $g(X):=\frac{(-K_{X})^3}{2}+1=12$.
These Fano threefolds were
quite mysterious objects and the attempt to 
find a geometrical description of them led Mukai to find 
their relationship with
the concept of varieties of
power sums.
First we recall the following definition:

\begin{defn}
Let $V$ be a $(v+1)$-dimensional vector space and 
let $F\in S^m \check{V}$ be a homogeneous forms of degree $m$ on $V$. Set
\[
\VSP(F,n)^o:=
\{(H_1,\dots, H_n)\mid H_1^m+\cdots+H_n^m=F\}
\subset \Hilb^n \mP(\check{V}).
\]
The closed subset $\VSP(F,n):=\overline{\VSP(F,n)^o}$ 
is called  the {\em{varieties of power sums}} of $F$.
\end{defn}

Mukai discovered the following beautiful description of prime Fano threefolds
of genus $12$
(\cite{Mu2}, \cite{Mukai12}. See also \cite{Schr}):
\begin{thm}\label{v22} Let $\{F_4=0\}\subset \mP(V)$ 
    be a general plane quartic curve. Then 
    \begin{enumerate}[$(1)$]
\item
 ${\rm{VSP}}(F_4,6)\subset \Hilb^{6} \mP(\check{V})$
is a general prime Fano threefold of genus $12$; and conversely,
\item every general prime Fano threefold of genus $12$ is of this form.
\end{enumerate}
\end{thm}

Again the main character is played by ineffective theta characteristics
because the Hilbert scheme of lines on $X$ is isomorphic to 
a smooth curve $\sH_{1}$ of genus $3$ and Mukai proved that 
the correspondence on $\sH_1\times \sH_1$ defined 
by intersections of lines on $X$ gives
an ineffective theta characteristic $\theta$ on $\mathcal{H}_1$.
More precisely,
$\theta$ is constructed so that
\[
I_{\theta}=\{(l,m)\in \sH_1\times \sH_1
\mid l\cap m\not =\emptyset,l\not =m\}.
\]
(See \ref{subsection:introtheta} for
more detailed explanations in our setting).
Now, by the result of Scorza recalled above,
there exists the Scorza quartic $\{F=0\}$
of the pair $(\sH_1,\theta)$
in the same ambient plane as the canonically embedded $\sH_1$. Mukai proved that
$X$ is recovered as $\VSP(F, 6)$.
This is the result (2) of Theorem \ref{v22}. The result (1) follows
from (2) since the number of the moduli of prime Fano threefolds
of genus $12$ is equal to
$\dim \sM_3=6$.

\subsection{Scorza quartics}
\label{subsection:IntroSc}~

Up to now the main idea to study spin curves of genus $g$
with ineffective theta
has been to try to associate to them a quartic hypersurface.
In the case $g=3$, this association turns out to be
the inverse of the Scorza map.

Here we would like to recall the results of 
Dolgachev and Kanev \cite[\S 9]{DK} for a
modern account of Scorza's
beautiful construction of this quartic hypersurface \cite{Sco3}. 

Following 
\cite[7.1.4 p.279]{DK} let $\Gamma \subset \mP^{g-1}$ be a canonical curve 
of genus $g$, $\theta$ an ineffective theta characteristic on it and
$I_{\theta}\subset \Gamma\times\Gamma$ the Scorza correspondence. 
We denote by $I_{\theta}(x)$ the fiber of $I_{\theta}\to \Gamma$ over $x$
and call it the {\em{theta polyhedron}} attached to $x$.
In other words, $I_{\theta}(x)$ is the unique member of
$|\theta +x|$. 
Since the linear hull $\langle I_{\theta}(x)-y\rangle$
  is a hyperplane of $\mP^{g-1}$, then we can define a morphism
  $\pi_{\theta}\colon I_{\theta}\rightarrow |\omega_{\Gamma}|=\check{\mP}^{g-1}$ 
as a composition of
  the natural embedding 
$I_{\theta}\hookrightarrow \Theta_{\Gamma}$ and the Gauss map
  $\gamma\colon\Theta^{{\rm {ns}}}_{\Gamma}\to\check{\mP}^{g-1}$, where
  $\Theta_{\Gamma}\subset J(\Gamma)$ is the theta divisor and 
$\Theta^{{\rm{ns}}}_{\Gamma}$ is the nonsingular locus of $\Theta_{\Gamma}$.
 Set-theoretically $\pi_{\theta}$ is the map $(x,y)\mapsto  
\langle I_{\theta}(x)-y\rangle$. 
The hyperplane $\langle I_{\theta}(x)-y\rangle$ is called the {\it{face}} 
  of $I_{\theta}(x)$ opposed to $y$.

The following is an important invariant of 
$(\Gamma,\theta)$:
\begin{defn}
 \label{introdefn:Gamma}  
 The image $\Gamma(\theta)$ of the above morphism 
 $\pi_{\theta}\colon I_{\theta}\to \check{\mP}^{g-1}$ 
(with reduced structure)
  is called the {\it{discriminant
 locus}} of the pair $(\Gamma,\theta)$.
 \vspace{0.10in}
\end{defn}

By Definition \ref{introdefn:Gamma}, we have the following diagram:

\begin{equation}
\label{iniezionscorza}
\xymatrix{ & I_{\theta}\subset \Gamma\times\Gamma \ar[dl]_{{ \pi_{\theta}}}
\ar[dr]^{p}\\
 \Gamma(\theta)\subset\check\mP^{g-1} &  & \Gamma \subset\mP^{g-1}. }
\end{equation}
\noindent
Dolgachev and Kanev point out that,
to construct the Scorza quartic,
the following three conditions
are needed, which Scorza overlooked
(see \cite[(9.1) (A1)--(A3)]{DK})
: 
\begin{enumerate}[({A}1)]
 \item
the degree of the map $I_{\theta}\to\Gamma(\theta)$ is two, 
namely,
$\langle I_{\theta}(x')-y'\rangle=\langle I_{\theta}(x)-y\rangle$ implies 
$(x',y')=(x,y)$ or $(y,x)$,
\item $\Gamma(\theta)$ is not contained in a quadric, and
\item
  $I_{\theta}$ is smooth (this condition is
modified in \cite[5.5.3]{topic}).
\end{enumerate} 
From now on in this subsection,
we assume these conditions.

We can define:
\[
  \overline{D}_H:=\pi_{\theta *}p^*(H\cap \Gamma)
\]
\noindent
as a divisor,
where $H$ is an hyperplane of $\mP^{g-1}$.

It is not difficult to see 
$\deg \Gamma(\theta)=g(g-1)$
by (A1) and (A3)
(see \cite[Corollary 7.1.7]{DK}).
Moreover,
by $\deg ( H\cap \Gamma)=2g-2$ and $\deg p=g$,
it holds
$\deg p^*(H\cap \Gamma)=2g(g-1)$.
By (A1), it is easy to see
$\deg \overline{D}_H=2g(g-1)$.   
Therefore we may expect
that $\overline{D}_H$ is a quadric section of $\Gamma(\theta)$.
This is true (\cite[Proposition 9.2]{DK}):   
\begin{prop}
\label{prop:quadric}
$\overline{D}_H$ is cut out by a quadric in 
$\check{\mP}^{g-1}$.
\end{prop}
To show this, we need the assumption (A2).

From here on, we give an explanation of the Scorza
quartic slightly different from that of \cite{DK}
(but essentially the same).
We define the correspondence:
\[
\sD:=\{(q_1,q_2)\mid q_1\in \overline{D}_{H_{q_2}}\}\subset \Gamma(\theta)
\times \Gamma(\theta),
\]
where $H_q$ is the hyperplane of $\mP^{g-1}$
corresponding to $q\in \check{\mP}^{g-1}$.
It is easy to see that $\sD$ is symmetric.
By Proposition \ref{prop:quadric},
we see 
that $\sD$ is the restriction of
a symmetric $(2,2)$ divisor $\sD'$ of 
$\check{\mP}^{g-1}\times \check{\mP}^{g-1}$.
Let $\{\check{F}_4=0\}$
be the quartic hypersurface
obtained by restricting $\sD'$ to the diagonal of
$\check{\mP}^{g-1}\times \check{\mP}^{g-1}$.
The Scorza quartic is the `dual' quartic in $\mP^{g-1}$
of $\{\check{F}_4=0\}$.

To explain this more precisely, we give a quick review of some
generality of the theory of polarity. 
Set $V: = H^{0}(\Gamma,\omega_{\Gamma}\check{)}$. Each homogeneous form  
$F\in S^{4} \check{V}$ defines a linear map:
\begin{eqnarray*}
\ap_{F}\colon S^2 V& \to & S^{2} \check{V}\\
G & \mapsto & 
P_G(F).
\end{eqnarray*}
called the {\em{apolarity map}} (cf. \cite[Definition 1.5]{DK}),
which is nothing but the linear extension of iterating polar maps. If
$\ap_{F}$ is an isomorphism $F$ is called {\it {non-degenerate}} and
then the inverse isomorphism is given by a $\check{F}\in S^{4} V$, that is  
${\ap_{F}}^{-1}=\ap_{\check{F}}$. The form $\check{F}\in S^{4} V$ is 
called {\it{the dual form}} of $F$ \cite[\S 2.3]{doldual}.

It turns out that the constructed $\{\check{F}_4=0\}$ is non-degenerate
and we can take the dual $\{F_4=0\}$, which is the Scorza quartic.

We explain one of the important properties of the Scorza quartic.
By the theory of polarity and the definition of $\check{F}_4$, 
the fiber of $\sD\to \Gamma(\theta)$ over a point $q\in \Gamma(\theta)$
is defined by the second polar $P_{H_q^2}(\check{F}_4)$ of
$\check{F}_4$. Moreover, by definition of $\Gamma(\theta)$,
it is easy to derive that 
$P_{H_q^2}(\check{F}_4)=a b$ for some
$a, b\in \Gamma$ such that $(a,b)\in I_{\theta}$,
where $a,b \in {\mP}^{g-1}$ is considered
as a linear form on $\check{\mP}^{g-1}$.
By definition of the dual, we have
$P_{a,b}(F_4)=H_q^2$.
By this property and (\ref{rangouno}),
the association of the Scorza quartic is the inverse of the Scorza map
in the case $g=3$.
 
It is expected that the Scorza quartic is useful for the study of a spin curve
but no deep properties of the Scorza quartic were known.
Firstly, its construction were not so explicit as the reader who has
follows us till here has certainly verified.
Secondly, 
Scorza's construction itself depends on three assumptions (A1)--(A3)
as above and it were unknown whether these conditions are
fulfilled for a general spin curve of genus $>3$.
Thus even the existence itself of the
Scorza quartic was conditional except for the genus $g=3$ case,
solved by Scorza himself. 
The sum of \cite{TZ} and of this work solve these problems for every
$g\geq 3$ (see Theorem \ref{thm:main3}).

\subsection{Special quartics and incidence correspondences}
\label{subsection:Special}~

In \cite{TZ}, among other results, we constructed certain special quartics.
It is important to notice that in the end 
the construction of the
special quartics in \cite{TZ} is almost straightforward 
even if it relies on
geometrical ideas requiring a bit of technical work.  
The byproduct is that we can show the special quartics
coincide with the Scorza
quartics of trigonal spin curves (see \ref{subsection:introtheta}). 
For ease of reading
this paper, we think it useful to review results of \cite{TZ}.


Let $B$ be the smooth quintic del Pezzo threefold, 
that is $B$ is a smooth projective threefold
such that $-K_{B}=2H$, where $H$ is the ample generator
of $\Pic B$ and $H^3=5$. It is well known that the linear system 
$|H|$ embeds $B$ into $\mP^6$.
Recall that this image of $B$ can be seen 
as $G(2,5)\cap \mP^{6}$, where $\mP^{6}\subset \mP^9$ 
is transversal to $G(2,5)$ which is the Grassmannian
of the $2$-dimensional vector subspaces of a $5$-dimensional vector
space considered embedded into $\mP^{9}$ (see \cite{Fu2}, 
\cite[Thm 4.2 (iii), the proof p.511-p.514]{I1}).
We started from a general smooth rational curve 
$C$ of degree $d$ on $B$, 
where $d$ is an arbitrary integer greater than or equal to $6$
(see \ref{subsection:Cd} for more detailed properties of $C$). 
Let $f\colon A\to B$ be the blow-up along $C$ and
$E_{C}$ the $f$-exceptional divisor.
We define:
\begin{defn}
    \label{linea}
    A connected curve $l\subset {A}$ 
    is called a {\em{line}} on ${A}$ if
$-K_{{A}} \cdot l=1$ and
${{E_{C}}}\cdot l=1$.
\end{defn}

We point out that since $-K_A=f^*(-K_{B})-E_{C}$ and $E_{C}\cdot
l=1$ then
$f(l)$ is a line on $B$ {\it{intersecting}} $C$.
The classification of lines on $A$ is simple:
\begin{prop}\label{lineeA}
A line $l$ on $A$ is one of the following curves on $A:$
\renewcommand{\labelenumi}{\textup{(\roman{enumi})}}
\begin{enumerate}
\item
the strict transform of a uni-secant line of $C$ on $B$,
or
\item
the union $l_{ij}=\beta'_{i}\cup \zeta_{ij}$ $(i=1,\ldots, s,
\,j=1,2)$, where $\beta'_i$ is a bi-secant line $\beta_i$ of $C$ and
$\zeta_{ij}$ is the fiber of $E_C$ over a point in $C\cap \beta_i$.
\end{enumerate}
In particular $l$ is reduced and $p_a(l)=0$.
\end{prop}

\begin{prop}
\label{prop:intro1}
The Hilbert scheme of lines on $A$ 
is a smooth trigonal curve
$\mathcal{H}_1$ of genus $d-2$.
\end{prop}
We remind the reader of that Mukai constructed his plane quartics from 
$(\sH_1,\theta)$, which is a data of intersections of lines.
Instead, to construct the special quartics, 
we need the notion of conics on $A$ and data of their intersections.
\begin{defn}
    \label{conica}
    We say that a connected and reduced curve $q\subset {A}$ 
    is a {\em{conic}} on ${A}$ if
    $-K_{{A}} \cdot q = 2$ and ${{E_{C}}}\cdot q=2$.
\end{defn}

We showed that
the Hilbert scheme of conics on $A$ is an irreducible surface 
and the normalization morphism is injective, 
namely, the normalization
$\sH_2$ parameterizes conics on $A$ in one to one way.  

Moreover we have the full description
of $\sH_2$ as follows (\cite[Theorem 4.2.15]{TZ}).
For this,
let 
$D_l\subset \sH_2$ be the locus 
parameterizing conics on $A$ 
which intersect a fixed line $l$ on $A$.

\begin{thm}\label{thm:H_2}
$\sH_2$ is smooth and is 
a so-called {\em{White surface}} 
obtained by blowing up 
$S^2 C\simeq \mP^2$ at $s:=\binom{d-2}{2}$ points.
The locus $D_l$ is a divisor
linearly equivalent to
$(d-3)h-\sum_{i=1}^s e_i$ on $\sH_2$,
where $h$ is the pull-back of a line, $e_i$ 
are the exceptional curves of $\eta\colon \sH_2\to \mP^2$,
and
$|D_l|$ embeds
$\sH_2$ into $\check{\mP}^{d-3}$.
$\sH_2\subset \check{\mP}^{d-3}$ is 
projectively Cohen-Macaulay,
equivalently,
$h^i(\check{\mP}^{d-3}, \sI_{\sH_2}(j))=0$
for $i=1,2$ and $j\in \mZ$,
where $\sI_{\sH_2}$ is the ideal sheaf of $\sH_2$ in
$\check{\mP}^{d-3}$. 
Moreover, $\sH_2$ is given by intersection of cubics.  
\end{thm}

Here we use the notation $\check{\mP}^{d-3}$
since the ambient projective space of $\sH_2$
and that of the canonical embedding of $\sH_1$
can be considered as reciprocally dual
(see \ref{subsection:duality}).
We write the ambient of $\sH_1$ by $\mP^{d-3}$ and
that of $\sH_2$ by $\check{\mP}^{d-3}$.

Finally set
$$
\sD_2:=\{(q_1,q_2)\in\sH_{2}\times\sH_{2}\mid q_1\cap q_2\neq \emptyset\}
$$
\noindent 
and denote by $D_q$ the fiber of $\sD_2\to \sH_2$ over a point $q$.
Then $D_q\sim 2D_l$ and it holds that $\sD_2\sim p_1^*D_q+p_2^*D_q$.
In particular 
since 
$\sH_2$ is not contained in a quadric,
it holds
$H^0(\sH_2\times \sH_2,\sD_2)\simeq 
H^0(\check{\mP}^{d-3}\times \check{\mP}^{d-3},
\sO(2,2))$.
Thus
$\sD_2$ is the restriction of a unique $(2,2)$-divisor $\sD'_2$ on 
$\check{\mP}^{d-3}\times \check{\mP}^{d-3}$.
Since $\sD'_2$ is symmetric, we may assume 
its equation $\widetilde{\sD}_2$ is also symmetric.
The restriction of $\widetilde{\sD}_2$ to the diagonal is 
a quartic hypersurface $\{\check{F}'_4=0\}$ in $\check{\mP}^{d-3}$.
We showed that $\check{F}'_4$ is non-degenerate.
Then the desired quartic is the unique quartic hypersurface $\{F'_4=0\}$ 
in $\mP^{d-3}$ dual to $\check{F}'_4$.

Notice that this construction is quite similar to
that of the Scorza quartic.
This similarity will be clear once we define a theta characteristic
on $\sH_1$ and clarify the relation of $\sH_1$ and $\sH_2$
(see \ref{subsection:introtheta}).

The following is the main result of \cite{TZ}, which is also a
generalization of (2)
of Theorem \ref{v22}:
\begin{thm}
\label{thm:main2}
Let $f\colon A\to B$ be the blow-up along $C$, and
let $\rho\colon \widetilde{A}\to A$ be the blow-up
of $A$ along the strict transforms $\beta'_i$
of $\binom{d-2}{2}$ bi-secant lines $\beta_i$ of $C$ on $B$.
Then there is an injection from $\widetilde{A}$ to 
$\VSP(F'_4,n)$,
where $n:=\binom{d-1}{2}$.
Moreover the image of $\widetilde A$ is uniquely determined by
$\sD_{2}$ and is an irreducible component of
\[
\VSP(F'_4, n;\sH_2):=
\overline{
\{(H_1,\dots, H_n)\mid H_i \in \sH_2\}}
\subset \VSP(F'_4,n).
\]
\end{thm}

\subsection{Existence of the Scorza quartics}
\label{subsection:introtheta}~

Now we expose the results of this paper.
Consider the Scorza correspondence $I_{\theta}$
for a curve $\Gamma$ of genus $g$ and an ineffective theta $\theta$.
First notice the following,
which can be easily seen
by the Riemann-Roch theorem and very standard arguments:
\begin{enumerate}[(a)]
\item
$\theta=I_{\theta}(x)-x$ is (of course) independent of $x\in \Gamma$,
\item
$h^{0}(\Gamma, \theta +x)=1$ for any $x\in \Gamma$,
\item
$I_{\theta}$ is disjoint from the diagonal,
\item
$I_{\theta}$ is symmetric, and
\item
$I_{\theta}$ is a $(g, g)$-correspondence.
\end{enumerate}    
By \cite[Lemma 7.2.1]{DK}, conversely,
for any reduced correspondence $I'$ satisfying the
above conditions, there exists a unique ineffective 
theta characteristic such that $I'=I_{\theta}$. 

Now for the curve $\sH_1$ parameterizing lines on $A$
(see Proposition \ref{prop:intro1}), 
we can introduce the incidence correspondence: 
 \begin{equation}\label{toosimple}
I:=\{(l,m)\mid l\neq m, l\cap m\not=\emptyset\}
\subset \sH_1\times \sH_1
\end{equation}
with reduced structure.
Actually we need a more sophisticated way to define $I$: see \ref{subsection:theta}.
In Proposition \ref{prop:sopratheta} we prove
$I$ satisfies the conditions (a)--(e) whence
there exists a unique ineffective 
theta characteristic such that $I=I_{\theta}$. 
This is a generalization of Mukai's result explained above.

In \ref{subsection:duality}, we observe 
that there is a natural duality between $\sH_1$ and the space $\sH_2$ of
conics on $A$. This gives us a very computable way to produce 
the discriminant loci $\Gamma(\theta)$ of $\theta$: see Proposition 
\ref{subsection:discr1}. In particular, we prove that $\Gamma(\theta)$
is contained in $\sH_2$.

By virtue of our explicit computation of the discriminant,
we prove in Proposition \ref{prop:DK}
that the pair $(\mathcal{H}_1,\theta)$ satisfies the conditions 
\cite[(9.1) (A1)--(A3)]{DK}, which guarantee the existence of 
the Scorza quartic
for the pair $(\mathcal{H}_1,\theta)$.
Then, by a standard deformation theoretic argument,
we can then verify that the conditions (A1)--(A3) hold
also for a general spin curve,
hence we answer affirmatively to the Dolgachev-Kanev Conjecture:

\begin{thm}[=Theorem \ref{thm:DK}]
\label{thm:main3}
The Scorza quartic exists for a general even spin curve.
\end{thm}

Moreover we can find explicitly the Scorza quartic for $(\sH_1,\theta)$. 
In fact, by definition, the Scorza quartic for $(\sH_1,\theta)$
lives in $\mP (H^0(\sH_1,K_{\sH_1}\check{)})$ 
but by line-conic duality as in \ref{subsection:duality} we can consider
it lives in 
$\mP^{d-3}$.  
In \ref{subsection:coincide}, we prove 
\begin{prop}[=Proposition $\ref{prop:quadrics}$]
\label{prop:coincide}
The special quartic $\{F'_4=0\}\subset \mP^{d-3}$ of
Theorem $\ref{thm:main2}$ coincides with
the Scorza quartic of $(\sH_1,\theta)$.
\end{prop}

\subsection{Moduli spaces of trigonal even spin curves}~

Our construction has an application 
to the description of the moduli space $\sS^+_{d-2}$
of trigonal even spin curves.

In \cite[2.3]{TZ}, we constructed inductively a smooth rational curve $C_d$ on 
$B$ of degree $d$ by smoothing the union of a smooth rational curve $C_{d-1}$ of degree $d-1$ and a general uni-secant line of it on $B$.
We inductively define $\sH^{d}_{B}$ as 
the union of the components of the Hilbert scheme whose general point
parameterizes a smooth rational curve of
degree $d$ on $B$
obtained as smoothings of 
the unions of a general smooth rational curves of degree $d-1$
belonging to $\sH^{d-1}_B$ and their general uni-secant lines. 
Indeed, by \cite[Proposition 2.5.2]{TZ}, $\sH_d^B$ is irreducible.

It is known that $\Aut B$ is isomorphic to
the automorphism group $\PGL_2$ of 
the complex projective line.
The $\PGL_2$-action on $B$ induces the $\PGL_2$-action on $\sH_d^B$.
In the section \ref{section:moduli},
we show that
if $d\geq 7$ (resp.~$d=6$), 
then $\sS^+_{d-2}$ (resp.~$\sS^+_4$ or its double cover)
birationally parameterizes $\PGL_2$-orbits in $\sH^B_d$
(see Theorem \ref{thm:mod}).
   
In the forthcoming paper \cite{TZb},
we show that $\sS^+_4$ is rational using this description.
Indeed, we construct a $\PGL_2$-equivariant birational map 
$\sH^B_6\dashrightarrow (\sH^B_1)^6/\mathfrak{S}_6 \simeq (\mP^2)^6/\mathfrak{S}_6$, where $\mathfrak{S}_6$ is the symmetric group of degree $6$
acting on $(\mP^2)^6$ as the permutation of the factors, 
thus $\sS^+_4$ or its double cover 
birationally parameterizes $\PGL_2$-orbits of unordered six points 
in $\mP^2$.
We show the rationality of $\sS^+_4$
relating it with the classically studied moduli space
of $\PGL_3$-orbits of unordered six points in $\mP^2$.


\begin{ackn}
We are thankful to Professor S. Mukai for 
valuable discussions and constant interest on this paper.
We received various useful comments from 
K. Takeuchi, A. Ohbuchi, S. Kondo,
to whom we are grateful.
The first author worked on this paper
partially when he was staying at the Johns Hopkins University under the
program of Japan-U.S. Mathematics Institute (JAMI)
in November 2005 and at the Max-Planck-Institut f\"ur Mathematik
from April, 2007 until March, 2008. 
The authors 
worked jointly during the first author's stay at the Universit\`a di Udine
on August 2005, and the Levico Terme conference on Algebraic Geometry in
Higher dimensions on June 2007.
The authors are thankful to all the above institutes
for the warm hospitality they received. 
\end{ackn}
%
\section{Preliminaries}
In this section, mainly we review our results in the previous paper \cite{TZ}
which we need in this paper.

\subsection{Lines on the quintic del Pezzo threefold}
\label{subsection:lines}~

Let $\pi\colon \mP\to \sH_{1}^{B}$ be the universal
family of lines on the quintic del Pezzo threefold $B$ and 
$\varphi\colon \mP\to B$ the natural projection. 
By \cite[Lemma 2.3 and Theorem I]{FuNa}, 
$\mathcal{H}_{1}^{B}$ is isomorphic to $\mP^2$
and $\varphi$ is a finite morphism of degree three. 
In particular the number of lines passing through a point 
is three counted with multiplicities.  

Denote by 
$M(C)$ the locus $\subset \mP^2$ of lines 
intersecting an irreducible curve $C$ on
$B$,
namely, $M(C):=\pi(\varphi^{-1}(C))$ with reduced structure. 
Since $\varphi$ is flat, $\varphi^{-1}(C)$ is purely one-dimensional.
If $\deg C\geq 2$, 
then $\varphi^{-1}(C)$ does not contain a fiber of 
$\pi$, thus $M(C)$ is a curve. 

A line $l$ on $B$ is called a {\em{special line}} if
$\sN_{l/B}\simeq 
\sO_{\mP^1}(-1)\oplus \sO_{\mP^1}(1).$
Note that, if $l$ is not a special line on $B$,
then 
$\sN_{l/B}=\sO_{l}\oplus\sO_{l}$.

\begin{prop} 
\label{prop:FN}
It holds$:$ 
\begin{enumerate}[$(1)$]
\item Special lines are parameterized by a conic $Q_2$ on $\sH^1_B$,
\item
if $l$ is a special line,
then $M(l)$ is the tangent line to $Q_2$ at $l$.
If $l$ is not a special line,
then $\varphi^{-1}(l)$ is the disjoint union of 
the fiber of $\pi$ corresponding to $l$,
and the smooth rational curve dominating a line on $\mP^2$.
In particular,
$M(l)$ is the disjoint union of a line and the point $l\in \sH^1_B$.
{\em{By abuse of notation}}, we denote by $M(l)$ the one-dimensional part of
$M(l)$ for any line $l$. 
Vice-versa,
any line in $\mathcal{H}_{1}^{B}$ is of the form $M(l)$ for some line
$l$, and
\item 
the locus swept by lines intersecting $l$ is a hyperplane section $T_{l}$
  of $B$ whose singular locus is $l$. For every point $b$ of $T_l\setminus l$,
  there exists exactly one line which belongs to $M(l)$ 
  and passes through $b$.
Moreover, if $l$ is not special, then the normalization of $T_l$ is $\mF_1$
and the inverse image of the singular locus is the negative section of $\mF_1$,
or, if $l$ is special, then 
the normalization of $T_l$ is $\mF_3$
and the inverse image of the singular locus is the union of 
the negative section and a fiber.
\end{enumerate}
\end{prop}
\begin{proof} See \cite[\S 2]{FuNa} and \cite[\S 1]{IlievB5}. 
\end{proof}

\subsection{Smooth rational curves on the del Pezzo threefold}
\label{subsection:Cd}~

In \cite[2.3]{TZ}, we constructed inductively a smooth rational curve $C_d$ on $B$ of degree $d$
by smoothing the union of a smooth rational curve $C_{d-1}$ of degree $d-1$
and a general uni-secant line of it on $B$.

We inductively define $\sH^{d}_{B}$ as the union of
the components of the Hilbert scheme whose general point
parameterizes a smooth rational curve of
degree $d$ on $B$
obtained as smoothings of 
the unions of a general smooth rational curves of degree $d-1$
belonging to $\sH^{d-1}_B$ and their general uni-secant lines. 
Indeed, by \cite[Proposition 2.5.2]{TZ}, $\sH_d^B$ is irreducible.

A general $C_d$ belonging to $\sH^{d}_{B}$
has the following several nice properties:
 
\begin{prop}
\label{prop:Cd}
\begin{enumerate}[$(1)$]
\item
$\sN_{C_{d}/B}\simeq\sO_{\mP^{1}}(d-1)\oplus\sO_{\mP^{1}}(d-1)$. 
In particular 
$h^1(\sN_{C_{d}/B})=0$ and $h^0(\sN_{C_{d}/B})=2d$,
\item
there exist no $k$-secant lines of $C_d$ on $B$ with $k\geq 3$,
\item there exist at most finitely many bi-secant lines of $C_{d}$ on
$B$, any of them intersects $C_{d}$ simply,
and they are mutually disjoint,
\item neither a bi-secant line nor
a line through the intersection point between a bi-secant line and $C_d$
is a special line, 
\item 
$M(C_d)$ intersects $Q_2$ simply,
\item
$M(C_d)$ is an irreducible curve of degree $d$
with only simple nodes
$($recall that
we abuse the notation by denoting the one-dimensional part of
$\pi(\varphi^{-1}(C_1))$ by
$M(C_1)$$)$, and
\item
by letting $l$ be a general line intersecting $C_d$ or any bi-secant line
of $C_d$,  
$M(C_d)\cup M(l)$ has only simple nodes as its singularities.
\end{enumerate}
\end{prop}
\begin{proof} See \cite[Propositions 2.3.2, 2.4.1, 2.4.2 and 2.4.4]{TZ}.
\end{proof}

\subsection{Lines and conics on certain blow-ups of the del Pezzo 
threefold}
\begin{conv}
We usually denote by $\overline{l}$ the image of a line $l$ on $A$.
\end{conv}

Here we give a more precise definition of $D_l$
defined before Theorem \ref{thm:H_2}.
Inside $\sH_2\times \sH_1$, 
we can define the incidence loci:
\[
\widehat{\sD}_1:=\{(q, l)\in \sH_2\times \sH_1
\mid q \cap l \not =\emptyset\}.
\]
Let ${\sD}_1 \subset \sH_2\times \sH_1$ be 
the divisorial part of $\widehat{\sD}_1$. 
Since $\sH_1$ is a smooth curve ${\sD}_1 \to \sH_1$ is flat.  
Let ${D}_l$ be the fiber of 
${\sD}_1\to \sH_1$ over $l\in \sH_1$. Clearly we can write
$D_{l}\hookrightarrow \sH_{2}$.

The following result contains the nontrivial 
result that for a general $l\in\sH_{1}$, $D_{l}$ parameterizes
conics which properly intersect $l$.

\begin{prop}\label{noncontenere} 
    For a general $l\in\sH_{1}$, $D_{l}$
    does not contain any point corresponding to the line pairs $l\cup m$
    with $m\in\sH_{1}$, and hence $D_l$ parameterizes all conics
    which properly intersect $l$.
\end{prop}
\begin{proof}See \cite[Corollary 4.2.17]{TZ}. 
\end{proof}

\section{Existence of the Scorza quartic}
 
In this section we will use the geometries of the trigonal curve 
$\sH_{1}$ (see Proposition \ref{prop:intro1}) 
and of the White surface $\sH_{2}$ (see
Theorem \ref{thm:H_2}), respectively
to give an affirmative answer to the conjecture of Dolgachev and Kanev
 \cite[Introduction p. 218]{DK} (see Theorem \ref{thm:DK}).

\subsection{Scorza correspondence}
 \label{subsection:theta}~

We need a more scheme theoretic definition of the correspondence 
$I$ given in \ref{toosimple}. 

For a general $C_{d}\in \sH^{d}_{B}$, set
$C:= C_{d}$.
There is a natural morphism $\sH_1\to \sH_1^B\simeq \mP^2$
mapping the class of a line $l$ on $A$ to
that of the image $\overline{l}$ of $l$ on $B$.
The image of $\sH_1$ on $\sH_1^B$ is nothing but 
$M:=M(C)$ defined in \ref{subsection:lines}, 
and $\sH_1\to M$ is the normalization.
By Proposition \ref{prop:Cd} (6), $M$ has only nodes as its singularities.
By Proposition \ref{lineeA}, singularities of $M$
correspond to bi-secant lines of $C$.
Since $p_a(M)=\frac{(d-1)(d-2)}{2}$ and $g(\sH_1)=d-2$, 
the number of nodes of $M$,
which is equal to the number of bi-secant lines of $C$,
is $s:=\frac{(d-2)(d-3)}{2}$.

\begin{lem}
\label{lem:general}
$h^{0}(\sH_1,(\pi_{|\sH_1})^*\sO_M(1))=3$.
\end{lem}

\begin{proof}

Let $h\colon S \to \sH_1^{B}\simeq \mP^2$ be the blow-up of
$\sH_1^{B}$ at
the $s=\binom{d-2}{2}$ nodes of $M$. 
Then $\sH_1\sim d \lambda -2\sum_{i=1}^{s}\varepsilon_i$,
where $\lambda$ is the pull-back of a general line and $\varepsilon_i$ 
are exceptional curves.
By the exact sequence
\[
0\to \sO_S(\lambda-\sH_1) \to \sO_S(\lambda)  \to 
\sO_{\sH_1}((\pi_{|\sH_1})^*\sO_M(1))\to 0
\]
together with
$h^0(\sO_S(\lambda))=3$ and 
$h^0(\sO_S(\lambda-\sH_1))=h^1(\sO_S(\lambda))=0$,
we see that
$h^{0}(\sH_1,(\pi_{|\sH_1})^*\sO_M(1))=3$ is equivalent to
$h^1(\sO_S(\lambda-\sH_1))=0$.
By the Riemann-Roch theorem, we have
$\chi(\sO_S(\lambda-\sH_1))=0$. 
Thus by $h^0(\sO_S(\lambda-\sH_1))=0$,
$h^1(\sO_S(\lambda-\sH_1))=0$ is equivalent to
$h^2(\sO_S(\lambda-\sH_1))=0$.
By the Serre duality,
$h^2(\sO_S(\lambda-\sH_1))=h^0(\sO_S((d-4)\lambda-\sum_{i=1}^{s}\varepsilon_i)$.
Thus
we have only to prove that
there exists no plane curve of degree $d-4$ through $s$ nodes of $M$.
We prove this fact by using the inductive construction of $C=C_d$.
As we mentioned in \ref{subsection:Cd},
$C_{d+1}$ is obtained as 
the smoothing of the union of $C_d$ 
and a general uni-secant line $\overline l$ of $C_d$.
From now on in the proof, we put the suffix $d$ to the object
depending on $d$.
For example, $s_d:=\binom{d-2}{2}$.
In case $d=1$, the assertion is obvious.
Assuming  
$h^0(\sO_{S_d}((d-4)\lambda_d-\sum_{i=1}^{s_d}\varepsilon_{i,d})=0$,
we prove 
$h^0(\sO_{S_{d+1}}((d-3)\lambda_{d+1}-\sum_{i=1}^{s_{d+1}}
\varepsilon_{i,d+1})=0$.
By a standard degeneration argument,
we have only to prove that
there exists no plane curve of degree $d-3$ through $s_{d+1}$ nodes of 
$M_d\cup M(\overline{l})$, where $s_d$ of $s_{d+1}$ 
nodes are those of $M_d$ and
the remaining $s_{d+1}-s_d=d-2$ nodes 
are $M_d\cap M(\overline{l})$ except the two points
corresponding
to the two other lines $\overline l'$, 
$\overline l''$ through $C_d\cap \overline{l}$. 
Assume that 
there exists a plane curve $G$ of degree $d-3$ through $s_{d+1}$
nodes of $M_d\cup M(\overline{l})$. Then 
$G\cap M(\overline{l})$ contains at least $d-2$ points. 
Since $\deg G=d-3$, this implies $M(\overline{l})\subset G$. Thus 
there exists a plane curve of degree $d-4$ through $s_d$ nodes of
$M_d$, a contradiction.  
\end{proof}

We denote by $\delta$ the $g^1_3$ on $\sH_1$ which defines 
$\varphi_{|\sH_1}\colon \sH_1\rightarrow C$.
Let $l$, $l'$ and $l''$ be three lines on $A$ such that 
$l+l'+l''\sim \delta$. Then $\overline{l}$, $\overline{l}'$
and $\overline{l}''$
are lines through one point of $C$.
Set 
\[
\theta:=(\pi_{|\sH_1})^*\sO_M(1)-\delta.
\]
Note that $\deg \theta=d-3$.
 Let $l$ be any line on $A$ and
    $l', l''$ lines 
    such that $l+l'+l''\sim \delta$.  
    By $\theta+l=\pi_{|\sH_1}^*\sO_M(1)-l'-l''$ and
    Lemma \ref{lem:general}, we have 
    $h^{0}(\sH_1, \sO_{\sH_1}(\theta +l))=1$.
Let $p_i\colon \sH_1\times \sH_1 \to \sH_1$ ($i=1,2$) be the two
projections and 
    $\Delta$ the diagonal of $\sH_1\times \sH_1$.
    Set $\sL:=\sO_{\sH_1\times \sH_1}({p_2}^*\theta+\Delta)$. 
    By $h^{0}(\sH_1, \sO_{\sH_1}(\theta +l))=1$ for any $l\in
\sH_1$,
    we see that $p_{1*}\sL$ is an invertible sheaf.
    Define an ideal sheaf $\sI$ by ${p_1}^*p_{1*}\sL=\sL\otimes \sI$. 
    $\sI$ is an invertible sheaf and let $I$ be the divisor defined
by $\sI$. We will denote by $I(l)$ the fiber of $I\to \sH_1$ over $l$.
By definition,
$I(l)$ consists of
the points in the support of
$|\theta +l|$.
Since $\pi_{|\sH_1}^*\sO_M(1)-l'-l''$,
they correspond to 
lines on $B$ intersecting both $C$ and $\overline{l}$
except $\overline{l'}$ and $\overline{l''}$.
The number of them is at most $d-3$.
By Proposition \ref{prop:Cd} (7),
the number is actually $d-2$. 
Thus the fiber of $I\to \sH_1$ over a general $l$ is reduced. 
Hence $I$ is reduced. 

Now we show the following generalization of Mukai's result
\cite[\S 4]{Mukai12} in our setting:

\begin{prop}\label{prop:sopratheta} 
The class of 
    $\theta$ is an ineffective theta characteristic and
    $I=I_{\theta}$.
\end{prop}

\begin{proof}
    By invoking  \cite[Lemma 7.2.1]{DK} and the definition of $I$,
    it suffices to prove the following:
\begin{enumerate}[(a)]
\item
    $h^{0}(\sH_1, \sO_{\sH_1}(\theta +l))=1$ for any $l\in \sH_1$,
\item
$I$ is reduced,
\item
$I$ is disjoint from the diagonal,
equivalently, $(l,m)\in I$ if and only if $l\not =m$,
\item
$I$ is symmetric, and
\item
$I$ is a $(g(\sH_1), g(\sH_1))$-correspondence.
\end{enumerate}    

 Let $l$ be any line on $A$ and
    $l', l''$ lines 
    such that $l+l'+l''\sim \delta$.  

We have proved (a) and (b) already in the above discussion.

We prove (c).
It is equivalent to show that
the support of $I(l)$ does not contain $l$.
By definition $\theta +l=
\pi_{|\sH_1}^*\sO_M(1)-l'-l''$.
If $\overline{l}$ is special, then it is uni-secant
by Proposition \ref{prop:Cd} (4), and
$M$ is not tangent to $Q_2$ at $\overline{l}$
by Proposition \ref{prop:Cd} (5).
Hence we are done.
If $\overline{l}$ is not special, then
$M(\overline{l})$ does not contain $\overline{l}$, thus we are done.

We prove (d).
Let $m$ be a line on $A$ such that
$m\in \sH^1$ is contained in the support of
$I(l)$. It suffices to prove that for a general $l$,
$l\in \sH^1$ is contained in the support of
$I(m)$.
For a general $l$, we may assume that 
$m\not =l'$ or $l''$. 
Then it is easy to verify this fact.

Finally we prove (e).
Since $I$ is symmetric and $\deg (\theta+l)=d-2=g(\sH_1)$,
the divisor is a $(g(\sH_1),g(\sH_1))$-correspondence.    
\end{proof}

\subsection{Line-conic duality}
\label{subsection:duality}~

We consider the embeddings of $\sH_1$ and 
$\sH_2$ into projective spaces 
by the canonical linear system and 
the linear system $|D_l|$, respectively.
Then we show the ambient projective spaces are 
reciprocally dual. 

\begin{lem}
The projection $\sD_{1}\rightarrow\sH_{2}$
is finite and flat.
\end{lem}
\begin{proof}
Since 
$\sD_1$ is a Cartier divisor in a smooth threefold $\sH_1\times \sH_2$,
$\sD_{1}$ is Cohen-Macaulay. 
Since $M=M(C)$ is irreducible,
no conic on $A$ intersects infinitely many lines on $A$.
Therefore
$\sD_1\to \sH_2$ is finite, hence
$\sD_{1}\rightarrow\sH_{2}$ is flat
since $\sH_{2}$ is smooth.
\end{proof}
Denote by $\widetilde{H}_q$ the fiber of
the projection $\sD_{1}\rightarrow\sH_{2}$
over $q$.
For a general $q$, lines intersecting $q$ are general.
Thus, by Proposition \ref{noncontenere},
$\widetilde{H}_q$ parameterizes
all the lines intersecting a general $q$.
\begin{lem}
For a general conic $q$,
$\widetilde{H}_q\in |\pi^*\sO_M(2)-2\delta|$, namely,  
$\widetilde{H}_q\sim 2\theta\sim K_{\sH_1}$.
\end{lem}

\begin{proof}
Since $q$ is general, the image $\overline{q}$ of $q$ 
is a bi-secant conic of $C$.
Let $\overline{l}_i$ and $\overline{m}_j$ ($i=1,2,3, j=1,2,3$) be
the lines on $B$ through each point of $C\cap \overline{q}$
respectively.
Denote by $l_i$ and $m_j$ the lines on $A$ corresponding to
$\overline{l}_i$ and $\overline{m}_j$.
Since $q$ is general, lines $l_i$ and $m_j$ are also general.
By definition of $\delta$,
we have $l_1+l_2+l_3\sim m_1+m_2+m_3\sim \delta$.
The lines on $A$ intersecting $q$
come from lines on $B$ intersecting $C$ and $\overline{q}$
except $\overline{l}_i$ and $\overline{m}_j$ ($i=1,2,3, j=1,2,3$).
Therefore 
$\widetilde{H}_q\in |\pi^*\sO_M(2)-2\delta|$.  
\end{proof} 
By the flatness of $\sD_1\to \sH_2$,
it holds $\widetilde{H}_q\sim K_{\sH_1}$ for any $q$.

By Theorem \ref{thm:H_2},
$D_l$ is a hyperplane section of $\sH_2\subset \check{\mP}^{d-3}$.
Thus 
the family 
\[
\xymatrix
{
\sD_1\ar[r]\ar[d] & \sH_2\times \sH_1\ar[dl]\\
\sH_1 &
}
\]
induces the morphism 
\begin{eqnarray*}
\sH_1& \to & \mP^{d-3}\\
{l} & \mapsto & {D_l}
\end{eqnarray*}  
by the universal property of the Hilbert scheme,
where $\mP^{d-3}$ is the dual projective space
of $\check{\mP}^{d-3}$.
Since $D_l\not =D_{l'}$ for general $l\not =l'$, 
$\sH_1\to \mP^{d-3}$ is birational.
We denote by $\{H_q=0\}$ the hyperplane in $\mP^{d-3}$
corresponding to the point $q\in \check{\mP}^{d-3}$.
Note that,
for $l\in \sH_1$ and $q\in \sH_2$, 
$D_l\in \{H_q=0\}$ 
if and only if 
$D_l (q)=0$ by definition of $H_q$.
Thus $\widetilde{H}_q=\{H_q=0\}$ for a general $q$. 
Consequently, $\sH_1\to \mP^{d-3}$ 
coincides with the canonical embedding
$\Phi_{|K_{\sH_1}|}\colon \sH_1\to \mP^{d-3}$
by $\widetilde{H}_q\sim K_{\sH_1}$.

\subsection{Discriminant locus}
 \label{subsection:discr1}
 ~

We consider
$\sH_1\subset {\mP}^{d-3}$ and $\sH_2\subset \check{\mP}^{d-3}$. 
For the pair $(\sH_1,\theta)$,
we can interpret $\Gamma(\theta)$
by the geometry of lines and conics on $A$ 
as follows:

\begin{prop}
 \label{prop:discr}
For the pair $(\mathcal{H}_1,\theta)$, 
the discriminant locus $\Gamma(\theta)$ is contained in $\sH_{2}$, and
the generic point of the curve $\Gamma(\theta)$ 
parameterizes line pairs on $A$.  
\end{prop}       

\begin{proof}
  
  Take a general point $(l_1,l_2)\in I$, equivalently,
  take two general intersecting lines $l_1$ and $l_2$.
  $l_1\cup l_2$ is a conic and
  the lines corresponding to the points of $I(l_1)-l_2$ 
  are lines intersecting $l_1$ except $l_2$.
  Thus by discussions in \ref{subsection:duality},
  the point in $\check{\mP}^{d-3}$ corresponding to the hyperplane  
  $\langle I(l_1)-l_2\rangle$ is nothing but
  $l_1\cup l_2\in \mathcal{H}_2$. This implies the assertion.
\end{proof}

\begin{prop} 
\label{prop:gamma}
The curve $\Gamma(\theta)$ belongs to 
the linear system $|3(d-2)h-4\sum_{i=1}^{s}e_i|$ on $\sH_2$.
In particular $\Gamma(\theta)$ is not contained in 
a cubic section of $\sH_2$.
\end{prop}

\begin{proof} 
For a point $b\in C$, set
\[
L_b:=\overline{\{q\in \sH_2 \mid \exists b'\not =b,
f(q)\cap C=\{b, b' \} \}}.
\]
We show that the image $\eta(L_b)$ of $L_b$ on $S^2 C$ is a line.
Choose $b'\in C$ 
such that there exists no line on $B$ through $b$ and $b'$.
By \cite[Corollary 3.2.3]{TZ},
there exists a unique conic on $B$ through $b$ and $b'$.
This implies that $\eta(L_b)$ is a line.

We can write:
\[
\Gamma(\theta) \sim ah-\sum m_i e_i,
\]
where $a\in \mZ$ and $m_i\in \mZ$.

For a general $b\in C$, $L_b$ intersects $\Gamma(\theta)$
simply.
Thus $a$ is the number of line pairs whose images on $B$ pass through $b$.
By noting there exists three lines $l_1$, $l_2$ and $l_3$ through $b$,
it suffices to count the number of reducible conics on $B$
having one of $l_i$ as a component except $l_1\cup l_2$,
$l_2\cup l_3$ and $l_3\cup l_1$. Thus $a=3(d-2)$.

Now We count the number
of line pairs belonging to $e_i$.
Each of such line pairs is 
of the form $l_{ij;k} \cup l_{ij}$,
where $l_{ij;k}$ ($k=1,2$) is 
the strict transform of the line through 
one of the two points in $\beta_i\cap C$
distinct from $\beta_i$.
Thus the number of such pairs is four whence $m_i\geq 4$.

Finally
we count the number
of line pairs intersecting
a general line $l$. By Proposition \ref{noncontenere}, $D_{l}$ does not
contain any line pair having $l$ as a component. 
Since the number of lines on $A$ intersecting a fix line
on $A$ is $d-2$, 
we see that $D_{l}\cdot\Gamma(\theta)\geq (d-2)(d-3)$. Then
\[
(d-2)(d-3)\leq \Gamma(\theta) \cdot D_l=
(d-3)a-\sum_{i=1}^{s} m_i.
\]
where $s=\frac{(d-2)(d-3)}{2}$. Since we have shown $m_i\geq 4$,
this implies $m_i=4$.

\end{proof}

\begin{cor}
\label{cor:DK0}
For $(\sH_1,\theta)$,
it holds that
$\deg \Gamma(\theta)=g(g-1)$ and
$p_a(\Gamma(\theta))=\frac{3}{2}g(g-1)+1$.
Moreover, $K_{\Gamma(\theta)}=\sO_{\Gamma(\theta)}(3)$.
\end{cor}
\begin{proof}
The invariants of $\Gamma(\theta)$ are easily calculated
by Proposition \ref{prop:gamma}.
\end{proof}

\begin{cor}
\label{cor:DK1}
The restriction map $H^0(\sO_{\check{\mP}^{d-3}}(2))\to H^0(\sO_{\Gamma(\theta)}(2))$
is an isomorphism.
\end{cor}

\begin{proof}
By Theorem \ref{thm:H_2},
$H^0(\sO_{\check{\mP}^{d-3}}(2))\to H^0(\sO_{\sH_2}(2))$
is an isomorphism. 
To see 
$H^0(\sO_{\sH_2}(2))\to H^0(\sO_{\Gamma(\theta)}(2))$
is an isomorphism,
we have only to show that $H^1(\sH_2, \sO_{\sH_2}(2)\otimes
\sO_{\sH_2}(-\Gamma(\theta)))=\{0\}$.
By the Serre duality, the last cohomology group is isomorphic to
$H^1 (\sH_2, \sO_{\sH_2}(-2)\otimes
\sO_{\sH_2}(K_{\sH_2}+\Gamma(\theta))$,
and moreover, 
by $K_{\sH_2}+\Gamma(\theta)=\sO_{\sH_2}(3)$,
it is isomorphic to
$H^1(\sH_2, \sO_{\sH_2}(1))$, which vanishes by Theorem \ref{thm:H_2}.
\end{proof} 


\subsection{Existence of the Scorza quartic}\-
\label{subsection:conj}~

We show the three conditions as in \ref{subsection:IntroSc}
hold for general pairs of canonical curves $\Gamma$ and
ineffective theta characteristics $\theta$ 
as Dolgachev and Kanev conjectured.

First we show that for our trigonal curve $\sH_1$ and 
the ineffective theta characteristic $\theta$ defined by intersecting lines 
on $A$ the above conditions hold.

\begin{prop}
\label{prop:DK}
$(\sH_1,\theta)$ satisfies $(\mathrm{A}1)$--$(\mathrm{A}3)$.
\end{prop}

\begin{proof}
(A1) This condition means that for general lines $l$ and $l'$ on $A$ 
such that $(l,l')\in I$ the face $\langle I(l)-l'\rangle$ belongs only
to $I(l)$ and to $I(l')$.

By contradiction assume that there exists
a line $m$ on $A$ such that $m\neq l$, $m\neq l'$ and 
  $\langle I(l)-l'\rangle$ is a face of $I(m)$.
  Then  some $d-3$ points of $I(m)$ lie on the hyperplane
  $\langle I(l)-l'\rangle$, equivalently,
  $m$ intersects $d-3$ lines on $A$ 
  corresponding to $d-3$ points of 
  $I(l)\cup I(l')$ except $l$ and $l'$.
By $d\geq 6$,
  it holds that, for $l$ or $l'$, say, $l$,  there exist
  two lines intersecting both $l$ and $m$.

  Consider the projection $B\dashrightarrow Q$ from
  the line $f(l)={\overline l_{}}$.
By \cite{Fu2},
  the target of the projection is 
  the smooth quadric threefold $Q$ and
  the projection is decomposed as follows: 
\begin{equation*}
\label{eq:linebis}
\xymatrix{ & B_{\overline l}\ar[dl]_{\pi_{1}}
\ar[dr]^{\pi_{2}}\\
 B &  & Q,}
\end{equation*}
where $\pi_1$ is the blow-up along $\overline{l}$.
Moreover, the image $E'_{\overline{l}}$ of the $\pi_1$-exceptional divisor 
$E_{\overline{l}}$ on $Q$ is a hyperplane section.
 
Now notice that, 
by generality of $l$, ${\overline l}\not =\overline{m}:=f(m)$ is
equivalent to have $l\not =m$.
Assume by contradiction
that ${\overline l}\cap \overline{m}\not =\emptyset$.
Then they span a plane $P$, which
contains two lines intersecting both ${\overline l}$ and 
$\overline{m}$.
This implies that $P\subset B$
since $B$ is the intersection of quadrics,
a contradiction.
Thus ${\overline l}\cap \overline{m}=\emptyset$
whence 
the strict transform $\overline{m}'$ of $\overline{m}$ on $Q$ is a line. 
Since there exist two lines intersecting both ${\overline l}$ and 
$\overline{m}$, 
$\overline{m}'$ intersects the image $E'_{\overline l}$ of 
$E_{\overline l}$ at two points.
Since $E'_{\overline l}$ is a hyperplane section on $Q$,
this implies that $\overline{m}'\subset E'_{\overline l}$, a contradiction.
\\
(A2) This condition is satisfied by
Theorem \ref{thm:H_2} and Proposition \ref{prop:gamma}.\\
(A3) By \cite[Lemma 7.1.3]{DK},
$(m_1,m_2)\in I$ is a singular point of $I$
if and only if 
$|I(m_1)-2m_2|\not =\emptyset$
and
$|I(m_2)-2m_1|\not =\emptyset$.

Let $m$ be a line on $A$, and $l_1$ and $l_2$ two lines on $A$
such that $\delta\sim m+l_1+l_2$.
By definition of $\theta$,
$I(m)\sim \theta+m\sim (\pi_{|\sH_1})^*\sO_M(1)-l_1-l_2$.
Therefore 
$|I(m)-2n|\not =\emptyset$
if and only if one of the following holds:
\begin{enumerate}[(1)]
\item $\overline{n}$ is a smooth point of $M$.
In this case, $\overline{n}$ is a uni-secant line of $C$.
If $\overline{n}\not =\overline{l}_1$ nor
$\overline{l}_2$, then
$M(\overline{m})$ is tangent to $M$
at $\overline{n}$.
If $\overline{n}=\overline{l}_1$ or
$\overline{l}_2$, then
$M(\overline{m})$ is tangent to $M$
at $\overline{n}$ with multiplicities three, or
\item
$\overline{n}$ is a singular point of
$M$, which is a node. 
In this case, $\overline{n}$ is a bi-secant line of $C$.
Correspondingly, there is another line $n'$ on $A$, see proposition
\ref{lineeA} $(ii)$.
The two branches of $M$ at $\overline{n}$
correspond to $n$ and $n'$ respectively
since $\sH_1\to M$ is the normalization.
If $\overline{n}\not =\overline{l}_1$ nor
$\overline{l}_2$, then
$M(\overline{m})$ is tangent at $\overline{n}$ 
to the branch of $M$ corresponding to $n$.
If $\overline{n}=\overline{l}_1$ or
$\overline{l}_2$, then
$M(\overline{m})$ is tangent at $\overline{n}$ 
to the branch of $M$ corresponding to $n$
with multiplicity three.
\end{enumerate}
Recall that,
for a line $\overline{l}$ on $B$,
we denote by $T_{\overline{l}}$ the hyperplane section
swept out by lines intersecting $\overline{l}$
(Proposition \ref{prop:FN} (3)).
We can restate the above conditions as follows:
\begin{enumerate}[(1)]
\item
If $\overline{n}\not =\overline{l}_1$ nor
$\overline{l}_2$, then
$C$ is tangent to $T_{\overline{m}}$ at
$C\cap \overline{n}$.
Assume that $\overline{n}=\overline{l}_1$ or
$\overline{l}_2$. 
If $\overline{n}$ is not a special line, then
$C$ is tangent at $C\cap \overline{n}$
with multiplicity three
to the branch of $T_{\overline{m}}$ 
corresponding to $\overline{n}$.
If $\overline{n}$ is a special line, then
$C$ intersects $T_{\overline{m}}$ at $C\cap \overline{n}$
with multiplicity three.
\item
Note that, by Proposition \ref{lineeA}, 
$n$ corresponds to one of a point $p_n$ of
$C\cap \overline{n}$.
By Proposition \ref{prop:Cd} (4),
$\overline{n}$ is not a special line.
If $\overline{n}\not =\overline{l}_1$ nor
$\overline{l}_2$, then
$C$ is tangent to $T_{\overline{m}}$ at $p_n$.
If $\overline{n}=\overline{l}_1$ or
$\overline{l}_2$, then
$C$ is tangent at $p_n$
with multiplicity three
to the branch of $T_{\overline{m}}$ 
corresponding to $\overline{n}$.
\end{enumerate}

Bearing this in mind,
we prove that $I$ is smooth for a general $C$
by simple dimension count.
We only prove
$I$ is smooth at $(\overline{m}_1,\overline{m}_2)$
with both $\overline{m}_1$ and $\overline{m}_2$ non-special.
The remaining cases can be treated similarly.
Let $\overline{m}_1$ and $\overline{m}_2$ be two intersecting non-special
lines on $B$.
We estimate the codimension in $\sH^B_d$ of the locus $\sH'$ of $C$
such that $C$ intersects both $\overline{m}_1$ and $\overline{m}_2$
and is tangent to both $T_{\overline{m}_1}$ and $T_{\overline{m}_2}$.  
By Proposition \ref{prop:Cd} (1),
passing through one point is a codimension two condition.
Moreover, being tangent to a smooth surface is a 
codimension one condition.
The choice of two points on $\overline{m}_1$ and $\overline{m}_2$
respectively has two parameters.
Thus $\codim \sH'=4$.
Since the choice of $\overline{m}_1$ and $\overline{m}_2$
has three parameters,
we have the claim for
a general $C$. 
\end{proof}

For any spin curve $(\Gamma,\theta)$ with ineffective $\theta$,
let 
\[
\Gamma'(\theta):=I_{\theta}/(\tau),
\]
where $\tau$ is the involution on $I_{\theta}$ induced by
that of $\Gamma\times \Gamma$ permuting the factors. 
Note that $I_{\theta}\to \Gamma(\theta)$ factor through $\Gamma'(\theta)$.

\begin{cor}
\label{cor:DK2}
For $(\sH_1,\theta)$,
it holds
$\Gamma'(\theta)\simeq \Gamma(\theta)$.
In particular, $\Gamma(\theta)$ is a smooth curve.
\end{cor}

\begin{proof}
By Proposition \ref{prop:DK}, 
(A1) and (A3) hold for $(\sH_1,\theta)$.  
Thus
we have 
$p_a(\Gamma'(\theta))=
\frac{3}{2}g(g-1)+1$
by
\cite[Corollary 7.1.7]{DK}.
Thus
$p_a(\Gamma'(\theta))= p_a(\Gamma(\theta))$
by Corollary \ref{cor:DK0}.
By (A1) again,
the natural morphism
$\Gamma'(\theta)\to \Gamma(\theta)$ is birational.
Therefore
it holds
$\Gamma'(\theta)\simeq \Gamma(\theta)$.

Since $I$ is smooth, and $I$ is disjoint from the diagonal,
the map $I\to \Gamma'(\theta)$ is \'etale.
Thus $\Gamma(\theta)\simeq \Gamma'(\theta)$ is a smooth curve.
\end{proof}

By a moduli theoretic argument we prove the conjecture for a general pair 
$(\Gamma,{\theta})$.
\begin{thm}
\label{thm:DK}
A general spin curve   
satisfies the conditions $(\textup{A}1)$--$(\textup{A}3)$.
In particular, the Scorza quartic exists for a general spin curve.
\end{thm}

\begin{proof}
Classically, the moduli space $\sS^{+}_{g}$ of even spin curves of genus $g$
is known to be irreducible (see \cite{ACGH}). 
Let $U$ be a suitable finite cover of an open neighborhood of a
general
$(\sH_1,\theta)\in \sS^{+}_{g}$
such that
there exists the family
$\sC\rightarrow U$ of 
pairs of canonical curves 
and ineffective theta characteristics.
Denote by $(\Gamma_u,\theta_u)$ the fiber of $\sC\to U$ over $u\in U$. 
By Proposition \ref{prop:DK}, 
$(\sH_1,\theta)$ satisfies (A1)--(A3).
Since the conditions (A1) and (A3) are open conditions,
these are true on $U$. 
Thus we have only 
to prove that the condition (A2) is still true on $U$.
Let $\sJ\rightarrow U$ be the family of Jacobians and
$\Theta\rightarrow U$ the corresponding family 
of theta divisors. By \cite[p.279-282]{DK}, the family
$\sI$ of the Scorza correspondences embeds into $\Theta$, and 
by the family of Gauss maps
$\Theta\rightarrow \check{\mP}^{g-1}\times U$, 
we can construct the family
${\sG}\rightarrow U$ whose fiber
${\sG}_{u}\subset \check{\mP}^{g-1}$ 
is the discriminant $\Gamma(\theta_u)$.
By Corollary \ref{cor:DK2},
it holds $\Gamma'(\theta)\simeq\Gamma(\theta)$ for $(\sH_1,\theta)$.
Thus we have also $\Gamma'(\theta_u)\simeq \Gamma(\theta_u)$ for $u\in U$.
By \cite[Corollary 7.1.7]{DK},
we see that
$p_a(\Gamma(\theta_u))$ and $\deg \Gamma(\theta_u)$ are constant for $u\in U$. 
Thus ${\sG}\rightarrow U$ is a flat family
since the Hilbert polynomials of fibers are constant.
Since no quadrics contain $\Gamma(\theta)$ for $(\sH_1,\theta)$,
neither does $\Gamma(\theta_u)$ for $u\in U$ 
by the upper semi-continuity theorem.
\end{proof} 

We have the following corollary to
the proof of Theorem \ref{thm:DK}:
\begin{cor}
Let $(\Gamma,\theta)$ be a general pair of 
a canonical curve $\Gamma$ 
and an ineffective theta characteristic $\theta$.
\begin{enumerate}[$(1)$]
\item
$\Gamma(\theta)$ is smooth.
\item
$\Gamma'(\theta)\simeq \Gamma(\theta)$.
\item
$K_{\Gamma(\theta)}=\sO_{\Gamma(\theta)}(3)$.
\item
The restriction morphism
$H^0(\sO_{\check{\mP}^{d-3}}(2))\to H^0(\sO_{\Gamma(\theta)}(2))$
is an isomorphism.
\end{enumerate}
\end{cor}

\begin{proof}
(1) follows from (A3) for $(\Gamma,\theta)$.
For the other, by the deformation theoretic argument in the proof of
Theorem \ref{thm:DK},
we have only to show the assertion for
a general $(\sH_1,\theta)$ constructed from
the incidence correspondence of lines on $A$.
This is true by Corollaries \ref{cor:DK0},
\ref{cor:DK1}, and \ref{cor:DK2}.
\end{proof}

\subsection{Scorza quartic of trigonal spin curves}
\label{subsection:coincide}
\begin{prop}
\label{prop:quadrics} 
The special quartic $F'_4$ as in $\ref{subsection:introtheta}$
is the Scorza quartic for $(\sH_1,\theta)$.
\end{prop}

\begin{proof} 
As in \ref{subsection:IntroSc}, 
the dual $\check{F}_4$ of the Scorza quartic is obtained
by restricting to the diagonal
the $(2,2)$ divisor on $\check{\mP}^{d-3}\times \check{\mP}^{d-3}$
coming from the correspondence
\[
\sD:=\{(q_1,q_2)\mid q_1\in \overline{D}_{H_{q_2}}\}
\subset \Gamma(\theta)
\times \Gamma(\theta).
\]
On the other hand,
the special quartic $\check{F}'_4$
are obtained by restricting to the diagonal
the $(2,2)$ divisor of $\check{\mP}^{d-3}\times \check{\mP}^{d-3}$
coming from the correspondence
\[
\sD_2:=\{(q_1,q_2)\mid q_1\cap q_2\neq \emptyset\}
\subset \sH_{2}\times\sH_{2}.
\]
Actually, $\check{F}'_4$ is determined by
the restriction of $\sD_2$ to $\Gamma(\theta)\times \Gamma(\theta)$
by Theorem \ref{thm:H_2} and Corollary \ref{cor:DK1}.
Therefore the assertion is equivalent to show 
$\overline{D}_{H_q}=\{\widetilde{D}_q=0\} \cap \Gamma(\theta)$
for a general $q$.
The set $\{\widetilde{D}_q=0\} \cap \Gamma(\theta)$ consists of points 
corresponding to the line pairs on $A$ intersecting $q$. 
By definition of $\overline{D}_{H_q}$,
it is rather straightforward to see  
the set $\overline{D}_{H_q}$ also consists of points 
corresponding to the line pairs
intersecting $q$. 
\end{proof}

\section{Moduli space of trigonal even spin curves}
\label{section:moduli}
Let $\sM^{\mathrm{tr}}_g$ and $\sS^{+\mathrm{tr}}_g$ be
the moduli space of trigonal curves of genus $g$ and
the moduli space of trigonal even spin curves of genus $g$,
respectively.
Denote by $\sH_d^{B}$ the Hilbert scheme of 
general smooth rational curves of degree $d$ on $B$
obtained inductively as in Proposition \ref{prop:Cd}.
By \cite[Proposition 2.5.2]{TZ},
$\sH_d^B$ is irreducible.
It is known that $\Aut B$ is isomorphic to
the automorphism group $\mathrm{PGL}_2$ of 
the complex projective line (see \cite{MU} and \cite{PV}).
The $\PGL_2$-action on $B$ induces the $\PGL_2$-action on $\sH_d^B$.
We have a natural rational map $\pi_{\sS}\colon  
\sH_d^B\dashrightarrow \sS^{+\mathrm{tr}}_{d-2}$
which maps a general $C_d$ to $(\sH_1,\theta)$
and is constant on general $\PGL_2$-orbits.
By taking suitable compactifications of $\sH^B_d$ and
$\sS^{+\mathrm{tr}}_{d-2}$,
a resolution of indeterminancy of $\pi_{\sS}$ and
the Stein factorization,
we have rational maps
$p_{\sS}\colon \sH^B_d\dashrightarrow \widetilde{\sS}^{+\mathrm{tr}}_{d-2}$
and $q_{\sS}\colon \widetilde{\sS}^{+\mathrm{tr}}_{d-2}\dashrightarrow
{\sS}^{+\mathrm{tr}}_{d-2}$ 
such that a general fiber of $p_{\sS}$ is connected and
$q_{\sS}$ is generically finite.
Then the $\PGL_2$-orbit of a general point of $\sH^B_d$
is contained in a fiber of $p_{\sS}$.
The purpose of this section is to show the following:

\begin{thm}
\label{thm:mod}
A general fiber of $p_{\sS}$ contains a $\PGL_2$-orbit
as an open dense subset.
If $d\geq 7$, then $q_{\sS}$ is birational onto the image.
If $d=6$, then the degree of $q_{\sS}$ is at most two.
\end{thm}

In other words,
if $d\geq 7$ (resp.~$d=6$), 
then $\sS^+_{d-2}$ (resp.~$\sS^+_{d-2}$ or its double cover)
birationally parameterizes $\PGL_2$-orbits in $\sH^B_d$.
   

Now we give three lemmas to prove Theorem \ref{thm:mod}.
First, as in Mukai's case (cf. the explanation about the proof of
Theorem \ref{v22} (1)),
we can reconstruct the threefold 
$\widetilde{A}$, whose definition is in the statement of Theorem \ref{thm:main2}, via the curve $\sH_1$ and
the ineffective theta characteristic $\theta$.
\begin{lem}
\label{prop:recovery}
The isomorphism class of 
$\widetilde{A}$ is recovered from 
the isomorphism class of $(\sH_1,\theta)$.
\end{lem}

\begin{proof}
From $(\sH_1,\theta)$, we can define
 $\Gamma(\theta)$ as in Definition \ref{introdefn:Gamma}
 and $F_4$ by Proposition \ref{prop:quadrics}.
 By Theorem \ref{thm:H_2} and Proposition \ref{prop:gamma},
 $\sH_2$ is recovered from $\Gamma(\theta)$ as
 the intersection of cubics containing $\Gamma(\theta)$.
 The divisor $\sD_2\subset \sH_2\times \sH_2$
 as in Theorem \ref{thm:main2}
 is recovered from the dual $\check{F}_4$.
 Thus, by Theorem \ref{thm:main2} and 
Proposition \ref{prop:quadrics},
 $\tA$ is recovered from $F_4$ and $\sH_2$.
 \end{proof}


To refine Lemma \ref{prop:recovery},
we show the following using the techinique of the Mori theory:

\begin{lem}
\label{prop:revise}
Suppose that $d\geq 7$.
Then, for a general $\tA$,
there exists a unique smooth rational curve $C$ of degree $d$ on $B$
up to the $\PGL_2$-action 
such that $\tA$ is obtained from $B$ by blowing up $C$ 
and the strict transforms of its bi-secant lines.

Suppose that $d=6$.
Then, for a general $\tA$,
there exists at most two smooth sextic rational curves 
$C^1, \dots, C^a$ $(a\leq 2)$
on $B$ up to the $\PGL_2$-action 
such that $\tA$ is obtained from $B$ by blowing up
one $C^i$ $(1\leq i\leq a)$ and 
the strict transforms of its bi-secant lines.
\end{lem}

\begin{rem}
In the forthcoming paper \cite{TZb},
we show that $a=2$ in case where $d=6$.
\end{rem}

\begin{proof}
Let $U$ be an open subset of $\sH^B_d$
such that, for $C\in U$,
\begin{itemize}
\item
$C$ is a smooth rational curve,
\item
$\sN_{C/B}\simeq \sO_{\mP^1}(d-1)\oplus \sO_{\mP^1}(d-1)$,
\item
all the bi-secant lines $\beta_i$ $(1\leq i\leq s)$ of $C$
are mutually disjoint, 
and
\item
$\sN_{\beta'_i/A}\simeq \sO_{\mP^1}(-1)\oplus
\sO_{\mP^1}(-1)$ 
for the strict transform $\beta'_i$ on $A$ of $\beta_i$
$(1\leq i\leq s)$.
\end{itemize}
It suffices to show the finiteness as above on $U$.
More precisely,
assume that
$\tA$ has two contractions
$\rho_j\colon \tA\to A_j$ ($j=1,2$) with the following properties:\\
(1)
$\rho_j$
contracts disjoint $s$ exceptional divisors
$E^j_i\simeq \mP^1\times \mP^1$ $(1\leq i\leq s$)
with $\sN_{\tA/E^j_i}\simeq \sO_{\mP^1\times \mP^1}(-1,-1)$
to rational curves,
and \\
(2) there exists a birational morphism $f_j\colon A_j\to B$
blowing down a $\mP^1$-bundle $E^j$ to a rational curve $C^j$ of degree $d$
with $\sN_{C^j/B}\simeq \sO_{\mP^1}(d-1)\oplus \sO_{\mP^1}(d-1)$.
In particular, $E^j\simeq \mP^1\times \mP^1$.

Then
we show that $d=6$, $\{E^1_i\}=\{E^2_i\}$ as sets
and $\rho_1$ and $\rho_2$ are the contractions
of $\{E^j_i\}$ along two different directions.

Since the conormal bundle of $E^j_i$ is ample,
there is an analytic contraction contracting
one $E^j_i$ only.
Therefore
it does not happen that $E^1_i\cap E^2_{i'}\not =\emptyset$ 
and $E^1_i\not = E^2_{i'}$
for some $i$, $i'$ since then 
$E^1_i\cap E^2_{i'}$ cannot be contracted on 
$E^1_i$ nor $E^2_{i'}$.
Assume by contradiction that
there is an $E^1_{i'}$ disjoint from all the $E^2_i$'s.
Let $\overline{E}^1_{i'}$ be the image of $E^1_{i'}$ on $A_2$.
It holds that $\overline{E}^1_{i'}\simeq \mP^1\times \mP^1$ and
$\sN_{A_2/\overline{E}^1_{i'}}
\simeq \sO_{\mP^1\times \mP^1}(-1,-1)$.
Hence 
there is an analytic contraction contracting
$\overline{E}^1_{i'}$ only.
Since $B$ does not contain a copy of $\mP^1\times \mP^1$,
$\overline{E}^1_{i'}\cap E^2$ is not empty.
This is a contradiction since
$E_2\simeq \mP^1\times \mP^1$ and
$\overline{E}^1_{i'}\cap E^2$ cannot be contracted on $E^2$.
Thus it holds $\{E^1_i\}=\{E^2_i\}$ as sets.
This implies that $A_1\dashrightarrow A_2$ is a flop.
Since $\rho(A_1)=\rho(A_2)=2$,
$\rho_1$ and $\rho_2$ are the contractions
of $\{E^j_i\}$ along two different directions.

Therefore we have the following diagram:
\begin{equation}
\label{eq:tworay}
\xymatrix{
& A_1 \ar[dl]_{f_1}  &\dashrightarrow & A_2 \ar[dr]^{f_2} & \\
 B  &  &  & & B.}
\end{equation}
For simplicity of the notation, we use 
the same notation for divisors on $A_1$ and
their strict transforms on $A_2$.

Let $H$ be the pull-back of the ample generator of $\Pic B$ by $f_1$
and $L$ the pull-back of the ample generator of $\Pic B$ by $f_2$.
Since $f_{2*} H$ is an effective divisor
on $B$ on the right hand side,
it holds that 
$f_{2*} H\sim pf_{2*} L$,
where $p$ is a positive integer.
Since $H$ is an $f_2$-divisor on $A_2$,
we can write 
\begin{equation}
\label{eq:num1}
H\sim pL-q{E}^2,
\end{equation}
where $q$ is a non-negative integer.
It holds that $-K_{A_1}=2H-E^1$ and similarly
$-K_{A_2}=2L-{E}^2$.
Thus
\begin{equation}
\label{eq:num2}
2H-E^1=2L-{E}^2. 
\end{equation}
By (\ref{eq:num1}) and (\ref{eq:num2}),
we have
\begin{equation}
\label{eq:num3}
E^1=(2p-2)L-(2q-1){E}^2.
\end{equation}
By (\ref{eq:num1}) and (\ref{eq:num3}),
it holds
\[
\text{$L=\frac{2q-1}{2q-p} H-\frac{q}{2q-p} E^1$, and 
${E}^2=\frac{2p-2}{2q-p} H-\frac{p}{2q-p} E^1.$}
\]
$l:=\frac{q}{2q-p}$ must be an integer,
thus we have $q=lq'$ and $p=q'(2l-1)$ with an integer $q'$.
Moreover
$\frac{2q-1}{2q-p}=\frac{2lq'-1}{q'}$ is also an integer,
we have $q'=1$.
Thus we have $p=2q-1$.
We compute 
$(-K_{A_1})^2 E^1$ in two ways:
first, 
$(-K_{A_1})^2 E^1=(2H-E^1)^2 E^1=-4H(E^1)^2+(E^1)^3$.
Since $f_1$ is the blow-up along a smooth rational curve of degree $d$,
we have $H(E^1)^2=-d$ and $(E^1)^3=-(2d-2)$.
Thus $(-K_{A_1})^2 E^1=2d+2$.
Second, note that $(-K_{A_1})^2 E^1=(-K_{A_2})^2 {E}^1$
since $A_1\dashrightarrow A_2$ is a flop
(see \cite[Lemma 3.1 (1)]{Taka} or
\cite[the proof of Corollary 9.3 (3)]{SShB2} for example).
By $-K_{A_2}=2L-E^2$ and 
${E}^1=4(q-1){L}-(2q-1) E^2$,
it holds that
$(-K_{A_2})^2 {E}^1
=(2{L}-E^2)^2(4(q-1){L}-(2q-1) E^2)$.
Since $f_2$ is the blow-up along a smooth rational curve of degree $d$,
we have ${L}(E^2)^2=-d$ and $(E^2)^3=-(2d-2)$ on $A_2$.
Thus
$(-K_{A_2})^2 {E}^1=
76q-78-(8q-6)d$.
Now we have the equality:
\[
76q-78-(8q-6)d=2d+2.
\]
Since $d=\frac{19q-20}{2q-1}$,
we obtain the following solution for $d\geq 6$:
\begin{itemize}
\item
$d=6$ and $q=2$,
\item 
$d=8$ and $q=4$, or
\item$d=9$ and $q=11$.
\end{itemize} 
We compute $(E^1)^3$ in two ways.
First, $(E^1)^3=-(2d-2)$ as mentioned above.
Second, since $A_1\dashrightarrow A_2$ flops
the strict transforms of bi-secant lines of $C^1$,
$(E^1)^3_{A_1}=({E}^1)^3_{A_2}+2^3 s$, where
$s:=\frac{(d-2)(d-3)}{2}$ is the number of bi-secant lines
(see for example \cite[Lemma 3.1 (1)]{Taka}, or
\cite[Proposition 3.5]{TZb}, where
we give a detailed proof).
It holds \[
({E}^1)^3_{A_2}=
(4(q-1){L}-(2q-1) E^2)^3_{A_2}=
320(q-1)^3-2(2q-1)^3-(8q-10)(2q-1)^2 d.
\]
Finally we have
\[
-(2d-2)=320(q-1)^3-2(2q-1)^3-(8q-10)(2q-1)^2 d+8s.
\]
It is easy to verify only $d=6$ and $q=2$ satisfies this equality.

\end{proof}

\begin{rem}
In \cite[Proposition 3.10]{TZb}, we show that
the diagram (\ref{eq:tworay}) as in
the proof of Lemma \ref{prop:revise} really exists
in the case where $d=6$.
\end{rem}

We refine Lemma \ref{prop:recovery}.

\begin{lem}
\label{lem:1to1}
There is one to one correspondence
between the isomorphism classes of 
$\widetilde{A}$ 
and the isomorphism classes of $(\sH_1,\theta)$.
\end{lem}

\begin{proof}
In any case, note that $(\sH_1,\theta)$ is determined from $(B,C)$.

If $d\geq 7$, then $(B,C)$ is recovered from $\tA$
by Lemma \ref{prop:revise} up to the $\PGL_2$-action,
thus $(\sH_1,\theta)$ is also recovered.
Therefore, by Lemma \ref{prop:recovery}, the assertion follows.

Suppose $d=6$.
As in the case where $d\geq 7$,
we have only to show that
$(\sH_1,\theta)$ is recovered from $\tA$.
By Lemma \ref{prop:revise} and its proof,
there are at most two choices of the contraction $\tA\to B$.
We use the notation of the proof of Lemma \ref{prop:revise}.
By $p=3$ and $q=2$, the equality (\ref{eq:num1}) is
$H\sim 3L-2{E}^2$, and 
the equality (\ref{eq:num3}) is
$E^1=4L-3{E}^2$.
From these equalities it is easy to see that
there is one to one correspondence
between the sets of lines on $A_1$ and lines on $A_2$.
Moreover, if two lines on $A_1$ intersect, then
the corresponding two lines on $A_2$ intersect, and vice versa.
Thus we can identify the Hilbert schemes of lines on $A_1$ and $A_2$,
and the theta characteristics on them.
Thus $(\sH_1,\theta)$ is recovered from $\tA$.
\end{proof}

\begin{proof}[The proof of Theorem $\ref{thm:mod}$]
We show that the isomorphism classes of $\widetilde{A}$ 
form an at least $(2d-3)$-dimensional family.
Indeed, 
an isomorphism of pairs
$(B,C^1)\to (B,C^2)$, where 
$C^1$ and $C^2$ are smooth rational curves of degree $d$ on $B$,
induces an isomorphism of the corresponding $3$-folds
$\tA_1$ and $\tA_2$ as in Theorem \ref{thm:main2}.
Conversely,
let $\widetilde{A}_1$ and $\widetilde{A}_2$
be two mutually isomorphic $3$-folds as in Theorem \ref{thm:main2}.
Choose a contraction $\tA_1\to B$ and 
a smooth rational curve $C^1$ of degree $d$ on $B$
as in Theorem \ref{thm:main2}.
Then, by an isomorphism $\iota\colon \widetilde{A}_1\to \widetilde{A}_2$,
we obtain the corresponding
contraction $\tA_2\to B$ and 
the smooth rational curve $C^2$ of degree $d$ on $B$.
Thus $\iota$ induces an isomorphism of pairs
$(B,C^1)\to (B,C^2)$. 
Since $\dim \sH_d^B=2d$ and 
$\dim \Aut B=\dim \PGL_2=3$,
the isomorphism classes of $(B,C)$,
where $C\in \sH^B_d$,
form an at least $(2d-3)$-dimensional family.
Therefore the isomorphism classes of $\widetilde{A}$ 
form a at least $(2d-3)$-dimensional family
since the choices of the contractions $\tA\to B$ are finite
by Lemma \ref{prop:revise}.
Now Lemma \ref{lem:1to1} implies that
the isomorphism classes of $(\sH_1,\theta)$ 
form a at least $(2d-3)$-dimensional family.
Since $\dim \sM^{\mathrm{tr}}_{d-2}=2d-3$
and a smooth curve has only a finite number of theta characteristics,
it holds that 
the isomorphism classes of $(\sH_1,\theta)$ 
form an open set of
an irreducible component of $\sS^+_{d-2}$
dominating $\sM^{\mathrm{tr}}_{d-2}$.
In particular, $\dim \Ima \pi_{\sS}=2d-3$.
Moreover,
a general $\PGL_2$-orbit in $\sH^B_d$ is $3$-dimensional
since the isomorphism classes of $(B,C)$,
where $C\in \sH^B_d$,
form a $(2d-3)$-dimensional family.
Therefore a general fiber of $p_{\sS}$ is $3$-dimensional
and
contains a general $\PGL_2$-orbit as a dense open subset.
Then $\widetilde{\sS}^{+\mathrm{tr}}_{d-2}$
birationally parameterizes $\PGL_2$-orbits in $\sH^B_d$.
Now the description of $q_{\sS}$ follows from
Lemmas \ref{prop:revise} and \ref{lem:1to1}.
\end{proof}

By the proof of Theorem \ref{thm:mod}, we have the following:

\begin{cor}
\label{cor:gen}
The image of $\pi_{\sS}\colon \sH^B_d\dashrightarrow \sS^{+tr}_{d-2}$ 
is an irreducible component of 
$\sS^{+tr}_{d-2}$ dominating $\sM^{\mathrm{tr}}_{d-2}$.
In particular a general $\sH_1$ is a general trigonal curve of genus $d-2$.
\end{cor}

\end{document}